\documentclass[graybox]{svmult}

\usepackage{mathptmx}       
\usepackage{helvet}         
\usepackage{courier}        
\usepackage{type1cm}        
\usepackage{makeidx}         
\usepackage{graphicx}        
\usepackage{multicol}        
\usepackage[bottom]{footmisc}
\usepackage{amsfonts}
\usepackage{amsmath}
\usepackage{amssymb}
\DeclareMathOperator{\re}{Re}
\DeclareMathOperator{\im}{Im}

\makeindex             
\def\C{{\mathbb C}}
\def\R{{\mathbb R}}
\def\N{{\mathbb N}}

\def\S{{\mathbb S}}

\def\ii{{\bf i}}
\def\la{\langle}
\def\ra{\rangle}

\def\ep{\epsilon}
\def\lg{L^2(\Gamma)}
\def\pa{\partial}

\def\nc{{\mathcal N}}
\def\oc{{\mathcal O}}
\def\curl{{\rm curl}\,}

\graphicspath{%
    {converted_graphics/}
    {/}
}
\begin{document}

\title*{Location and Weyl formula for the eigenvalues of some non
self-adjoint operators  }
\author{Vesselin Petkov}
\institute{Vesselin Petkov \at Institut de Math\'ematiques de Bordeaux, 351,
Cours de la Lib\'eration, 33405  Talence, France , \email{ petkov@math.u-bordeaux.fr}}
%
%
\maketitle

\abstract*{Each chapter should be preceded by an abstract (10--15 lines long) that summarizes the content. The abstract will appear \textit{online} at \url{www.SpringerLink.com} and be available with unrestricted access. This allows unregistered users to read the abstract as a teaser for the complete chapter. As a general rule the abstracts will not appear in the printed version of your book unless it is the style of your particular book or that of the series to which your book belongs.
Please use the 'starred' version of the new Springer \texttt{abstract} command for typesetting the text of the online abstracts (cf. source file of this chapter template \texttt{abstract}) and include them with the source files of your manuscript. Use the plain \texttt{abstract} command if the abstract is also to appear in the printed version of the book.}
\abstract{ We present a survey of some recent results concerning the location and the Weyl formula for the complex eigenvalues of two non self-adjoint operators. We study the eigenvalues of the generator $G$ of the contraction semigroup $e^{tG}, \: t \geq 0,$ related to the wave equation in an unbounded domain $\Omega$ with dissipative boundary conditions on $\partial \Omega$. Also one examines the interior transmission eigenvalues (ITE) in a bounded domain $K$ obtaining a Weyl formula with remainder for the counting function $N(r)$ of complex (ITE). The analysis is based on a semi-classical approach.}

\section{Introduction}
\label{sec:1}
\renewcommand{\theequation}{\arabic{section}.\arabic{equation}}
\setcounter{equation}{0}

 Let $P(x, D_x)$ be a second order  differential  operator with $C^{\infty}(K)$  real-valued coefficients in a bounded domain $K \subset \R^d, \: d \geq 2,$ with $C^{\infty}$ boundary $\partial K$. Consider a boundary problem
\begin{equation} \label{eq:1.1}
\begin{cases} P(x, D_x)u  = f \: {\rm in}\: K,\\
B(x, D_x) u = g \: {\rm on}\: \partial K, \end{cases}
\end{equation}
where $B(x, D_x)$ is a differential operator with order less or equal to 1 and the principal symbol $P(x, \xi)$ of $P(x, D_x)$ satisfies $p(x, \xi) \geq c_0 |\xi|^2,\: c_0 > 0.$
Assume that there exists  $0 < \varphi < \pi$ such that the problem
\begin{equation} \label{eq:1.2}
\begin{cases} (P(x, D_x) - z)u  = f \: {\rm in}\: K,\\
B(x, D_x) u = g \: {\rm on}\: \partial K. \end{cases}
\end{equation}
is parameter-elliptic for every  $z \in \Gamma_{\psi} = \{ z: \arg z = \psi \},\: 0 < |\psi|  \leq \varphi.$ Then following a classical result of Agranovich-Vishik \cite{AV}
we can find a closed operator $A$ with domain $D(A) \subset H^2(K)$ related to the problem (\ref{eq:1.1}). Moreover, for every closed angle $Q =\{z \in \C: \alpha \leq \arg z \leq \beta \}\subset \{z\in \C: |\arg z| < \varphi\}$ which does not contain $\R^+$ there exists $a_Q > 0$ such that  the resolvent $(A - z)^{-1}$ exists for $z \in Q, |z| \geq a_Q$ and the operator $A$ has a {\it discrete spectrum} in $\C$ with eigenvalues with finite multiplicities. 

Let $\{\lambda_j\}_{j=1}^{\infty}$ be the eigenvalues of $A$ ordered as follows
$$0 \leq |\lambda_1| \leq ...\leq |\lambda_m| \leq ....$$
In  general $A$ is not a self-adjoint operator and the 
analysis of the asymptotics of the counting function 
$$N(r) = \#\{ |\lambda_j| \leq r\}\: {\rm as}\: r \to + \infty,$$ 
where every eigenvalues is counted with its multiplicity, 
is a difficult problem. In particular, it is quite complicated to obtain a Weyl formula for $N(r)$ with a remainder and many authors obtained results which yield only the leading term of the asymptotics. On the other hand, even for parameter-elliptic boundary problems 
the result in \cite{AV} says that in any domain $0 < \psi < |\arg z| < \varphi $ we can have only finite number eigenvalues but we could have a bigger eigenvalues-free domains. To obtain a better remainder in the Weyl formula for $N(r)$ we must obtain a eigenvalues-free region outside some {\it parabolic neighborhood} of the real axis.

On the other hand, in mathematical physics there are many problems which are not parameter-elliptic. Therefore, the results of \cite{AV} cannot be applied and the analysis of the eigenvalues-free regions must be studied by another approach.

For the spectrum of non self-adjoint operators we have three important problems:\\

(I) Prove the discreteness of the spectrum of $A$ in some subset $U \subset \C$,\\

(II) Find eigenvalues-free domains in $\C$ having the form

$$|\im\: z | \geq ßC_{\pm\delta} (|\re\: z| +1)^{\delta_{\pm}}, \:\pm\re\: z \geq 0,\: 0 < \delta_{\pm} < 1,$$

(III)  Establish a Weyl asymptotic with remainder for the counting function
$$N(r) = c r^{d} + {\mathcal O}(r^{d - \kappa}),\: r \to \infty,\:0 < \kappa < 1.$$

In this survey we discuss mainly the problems (II) and (III) for two non self-adjoint operators related to the scattering theory. The problem (I) is easer to deal with and the analysis of (II) in many cases implies  that $A - z$ is a Fredholm operator for $z$ in a suitable region. We apply a new semi-classical approach for both problems (II) and (III). The analysis of (II) is reduced to the inversibility of a $h$-pseudo-differential operator, while for the asymptotic of $N(r)$ one exploits in a crucial way the existence of  parabolic neighborhood of the real axis containing the (ITE). The purpose of this survey is to present the recent results in \cite{V1}, \cite{V2}, \cite{PV}, \cite{PV1}, \cite{P1}, \cite{CPR2}, where the above problems are investigated by the same approach. We expect that our arguments can be applied to more general non self-adjoint operators covering the case of parameter-elliptic boundary problems (\ref{eq:1.2}).

\section{Two spectral problems related to the scattering theory}
\renewcommand{\theequation}{\arabic{section}.\arabic{equation}}
\setcounter{equation}{0}

{\bf I.} Let $K \subset \R^d,$ $d \geq 2$, be a bounded non-empty domain and let $\Omega = \R^d \setminus \bar{K}$ be connected. We suppose that the boundary $\Gamma$ of $\Omega$ is $C^{\infty}.$ 
Consider the boundary problem
\begin{equation} \label{eq:2.1}
\begin{cases} u_{tt} - \Delta_x u = 0 \,\: {\rm in}\: \R_t^+ \times \Omega,\\
\partial_{\nu}u - \gamma(x) u_t = 0\, \: {\rm on} \: \R_t^+ \times \Gamma,\\
u(0, x) = f_1, \: u_t(0, x) = f_2 \end{cases}
\end{equation}
with initial data $f = (f_1, f_2) \in H^1(\Omega) \times L^2(\Omega) = {\mathcal H}.$
Here $\nu(x)$ is the unit outward normal to $x \in \Gamma$ pointing into $\Omega$ and $\gamma(x)\geq 0$ is a $C^{\infty}$ function on $\Gamma.$ The solution of (\ref{eq:2.1}) is given by 
$$(u, u_t) = V(t)f = e^{tG} f,\: t \geq 0,$$
 where $V(t)$ is a contraction semi-group  in ${\mathcal H}$ whose  generator
$$ G = \Bigl(\begin{matrix} 0 & 1\\ \Delta & 0 \end{matrix} \Bigr)$$
has a domain $D(G)$ which is the closure in the graph norm of functions $(f_1, f_2) \in C_{(0)}^{\infty} (\R^n) \times C_{(0)}^{\infty} (\R^n)$ satisfying the boundary condition $\partial_{\nu} f_1 - \gamma f_2 = 0$ on $\Gamma.$ The spectrum of $G$ in $\re\: z < 0$ is formed by isolated eigenvalues with finite multiplicity (see \cite{LP2} for $d$ odd and \cite{P1} for $d$ even), while the continuous spectrum of $G$ coincides with $\ii \R.$
Next, if $Gf =\lambda f$ with $f = (f_1, f_2) \neq 0$, we get 

\begin{equation} \label{eq:2.2}
\begin{cases} (\Delta - \lambda^2) f_1 = 0 \:{\rm in}\,\, \Omega,\\
\partial_{\nu} f_1 -  \lambda \gamma f_1 = 0\: {\rm on}\,\, \Gamma. \end{cases}
\end{equation}

Thus if $\re \lambda < 0, \: f \neq 0,\:(u(t, x), u_t(t, x)) = V(t) f = e^{\lambda t} f(x) $, then $u(t, x)$ is a solution of (\ref{eq:2.1}) with {\it exponentially decreasing global energy.} Such solutions are called asymptotically disappearing and they perturb the inverse scattering problems. Recently it was proved (see \cite{CPR1}) that if we have a least one eigenvalue $
\lambda$ of $G$ with $\re \lambda < 0$, then the wave operators $W_{\pm}$ related to the problem (\ref{eq:2.1}) and the Cauchy problem for the wave equation are not complete, that is ${\text Ran}\: W_{-} \not=  {\text Ran}\: W_{+}$. Hence  we cannot define the scattering operator $S$ related to (\ref{eq:2.1}) by $S = W_{+} ^{-1}\circ W_{-}.$ We may define $S$ by another evolution operator. For problems associated to unitary groups, the associated scattering operator $S(z): L^2(\S^{d-1}) \to L^2(\S^{d-1})$ satisfies the equality
$$S^{-1}(z)= S^*(\bar{z}),\: z \in \C,$$
provided that $S(z)$ is invertible at $z$. This implies that $S(z)$ is invertible for $\im\: z > 0$, since $S(z)$ and $S^*(z)$ are analytic for $\im\: z < 0$ (see \cite{LP1} for more details). For dissipative boundary problems the above relation is not true and $S(z_0)$ may have a non trivial kernel for some $z_0, \im\: z_0 > 0.$ In the case of odd dimensions $d$ Lax and Phillips \cite{LP2} proved that $\ii z_0$ is an eigenvalue of $G$. Consequently, the analysis of the location of the eigenvalues of $G$ is important for the inverse scattering problems.
 
The eigenvalues of $G$ are symmetric with respect to the real axis, so it is sufficient to examine the location of the eigenvalues whose imaginary part is  nonnegative.
 A. Majda \cite{Ma} proved that
 if $\sup_{x \in \Gamma} \gamma(x) < 1,$ then the eigenvalues of $G$ lie in the region
$$E_1 = \{ z \in \C: \: |\re\: z | \leq C_1 (|\im\: z|^{3/4} + 1),\: \re\: z < 0\},$$
while if $\sup_{x \in \Gamma} \gamma(x) \geq 1,$ the eigenvalues of $G$ lie in $E_1 \cup E_2$, where
$$E_2 = \{ z \in \C: \: |\im\: z| \leq C_2 (|\re\: z|^{1/2} + 1),\: \re\: z < 0\}.$$
The case $\gamma(x) = 1,\: \forall x \in \Gamma,$ is special since as it was mentioned by Majda \cite{Ma} for some obstacles there are no eigenvalues of $G$. On the other hand, to our best knowledge we did not found a proof of this result in the literature. In the Appendix in \cite{P1}, the case when $K =  B_3 = \{x \in \R^3:\: |x| \leq 1\}$ is ball and $\gamma > 0$ is a constant has been examined and it  was proved that if $\gamma = 1$, there are no eigenvalues of $G$. On the other hand, for $\gamma = const > 1$ all eigenvalues of $G$ are real and for $0 < \gamma < 1$ there are no real eigenvalues. 

We will improve the above result of Majda and one examines two cases:
$$(A): \: 0 < \gamma(x) < 1,\: \forall x \in \Gamma.$$ 
 $$(B): \: \gamma(x) > 1, \:\forall x \in \Gamma.$$

{\bf II.} We discuss another important spectral problem for inverse scattering leading to non self-adjoint operator. For simplicity we assume that $d$ is odd. The inhomogeneous medium in $K$ is characterized by a smooth function $ n(x) > 0$ in $\bar{K}$, called {\it contrast}.  The scattering problem is related to an {\it incident wave} $u_i$ which satisfies the equation $(\Delta + k^2) u_i = 0$ in $\R^d$ and the 
{\it total wave}
$u = u_i + u_s$ satisfies the transmission problem
\begin{equation} \label{eq:2.3} 
\begin{cases} \Delta u + k^2 u = 0 \:\: {\rm in}\:\: \R^d \setminus \bar{K},\\
\Delta u + k^2 n(x) u = 0\:\: {\rm in}\:\: K,\\
u^{+} = u^{-} \:{\rm on}\: \Gamma,\\
\Bigl(\frac{\partial u}{\partial \nu}\Bigr)^{+} = \Bigl(\frac{\partial u}{\partial \nu}\Bigr)^{-}\:{\rm on}\: \Gamma, \end{cases},
\end{equation}
where $f^{\pm}(x) = \lim_{\epsilon \to 0}  f(x \pm \epsilon \nu(x))$ for $ x \in \Gamma$. Here $k > 0$  and the outgoing scattering wave $u_s$ satisfies the outgoing Sommerfeld radiation condition
$$\lim_{r \to + \infty} r^{(1-d)/2} \Bigl( \frac{\partial u_s}{\partial r} - \ii k u_s\Bigr)= 0$$
 uniformly with respect to $\theta = x/r \in \S^{d-1}, \: r = |x|$.

If the incident wave has the form $u_i = e^{\ii k \langle x, \omega \rangle},\: \omega \in \S^{d-1},$ then
$$u_s(r\theta, k)  = e^{\ii k r}r^{-(d-1)/2} \Bigl( a(k, \theta, \omega) + {\mathcal O} (\frac{1}{r})\Bigr), \: r \to +\infty.$$
The function $a(k, \theta, \omega)$ is called scattering amplitude and the {\it far-field operator} $F(k): L^2(\S^{d-1}) \longrightarrow L^2(\S^{d-1})$ has the form
$$(F(k)f)(\theta) = \int_{\S^{d-1}} a(k, \theta, \omega) f(\omega) d \omega.$$
 Notice also that the scattering operator has the
representation 
$$S(k)= I + \Bigl(\frac{\ii k}{2\pi}\Bigr)^{(d-1)/2} F(k).$$

The inverse scattering problem of the reconstruction of $K$  based on the linear sampling method  of Colton and Kress (see \cite{CH} ) breaks down for frequencies $k$ such that $F(k)$ has a non trivial kernel or co-kernel. 
Assume that for some $k \in \R^+$ the kernel of $F(k)$ is not trivial and let $F(k)f = 0,\: f \neq 0$.  We may consider an incident Herglotz wave 
$$u_i(x) = \int_{\S^{d-1}} e^{\ii k \langle x, \omega \rangle} f(\omega) d \omega.$$
Then one obtains
a scattering wave $u_s = {\mathcal O}(\frac{1}{r^2})$ since the  leading term 
$$\int_{\S^{d-1}} a(k, \theta, \omega) f(\omega) d\omega = 0$$
 vanishes. On the other hand, $(\Delta + k^2) u_s = 0$ in $\R^d \setminus \bar{K}$, so the Rellich theorem implies $u^s = 0$ in $\R^d \setminus \bar{K}$.
 Therefore the functions
$u = u_i\vert_{K} \neq 0$ and $w = (u_i +u_s)\vert_{K}$ satisfy the following problem 
\begin{equation}\begin{cases} \label{eq:2.4}
\Delta u + k^2 u = 0 \: {\rm in}\: K,\\
\Delta w + k^2 n(x) w = 0\: {\rm in}\: K,\\
u = w,\: \partial_{\nu} u = \partial_{\nu} w  \: {\rm on} \:\Gamma \end{cases}
\end{equation}
and $\lambda = k^2$ is called {\it interior transmission eigenvalue} (ITE). The inverse statement in general is not true and we may have complex (ITE).

 We  consider a more general setting. For $d \geq 2$, a complex number $\lambda\in \C \setminus \{0\},$ is called interior  transmission eigenvalue (ITE) if the following problem has a non-trivial solution $(u_1, u_2) \neq 0$:
\begin{equation}\label{eq:2.5}
\begin{cases} \left(\nabla c_1(x)\nabla+\lambda n_1(x)\right)u_1=0\: {\rm in}\:K,\\
\left(\nabla c_2(x)\nabla+\lambda n_2(x)\right)u_2=0\: {\rm in}\:K,\\
u_1=u_2,\,\,\, c_1\partial_\nu u_1=c_2\partial_\nu u_2\: {\rm on}\: \Gamma,\end{cases}
\end{equation}
where $c_j(x),n_j(x)\in C^\infty(\overline{K})$, $j=1,2$ are strictly positive real-valued functions.
For the analysis of (ITE) one imposes the condition
\begin{equation} \label{eq:2.6}
d(x) = c_1(x)n_1(x) - c_2(x)n_2(x) \neq 0,\quad\forall x\in\Gamma.
\end{equation}
Partial cases: 1) isotropic case: $c_1(x) = c_2(x),\: \forall x \in \Gamma$, $n_1(x) = 1, n_2(x) \neq 1, \: \forall x \in \Gamma.$ 2) anisotropic case: $c_1(x) \neq c_2(x),\: \forall x \in \Gamma.$ 

\section{Dirichlet-to-Neumann map}
\setcounter{equation}{0}

 The analysis of the eigenvalues-free domains is  based on a  semi-classical analysis.  Let $0 < h \ll 1$ and let $P(h) = -h^2 \Delta$. Introduce the sets

$$Z_1 = \{z \in \C:\: \re\: z = 1,  h^{1/2 - \epsilon} \leq \im\: z \leq 1,\: 0 < \epsilon \ll 1\},$$
$$Z_2 = \{z\in \C:\: \re\:z = -1, |\im\: z| \leq 1\},$$
$$Z_3 = \{z \in \C:\: |\re\: z | \leq 1, \im\: z = 1\}.$$
and consider for $z \in Z_1 \cup Z_2 \cup Z_3$ the semi-classical problem
\begin{equation} \label{eq:3.1}
\begin{cases} (P(h) - z)u = 0 \:\: {\rm in}\:\Omega,\: u \in H^2(\Omega),\\
u = f \:{\rm on}\: \Gamma, \end{cases}
\end{equation}

We need to introduce some $h$-pseudo-differential operators on a manifold with boundary $V$. We say that $a(x, \xi; h) \in S_{\delta}^{k, m}(V)$ if the following conditions are satisfied:

(i) for $|\xi| \geq L \gg 1$ we have\\
$$|\partial_x^{\alpha} \partial_{\xi}^{\gamma} a(x, \xi; h)| \leq C_{\alpha, \gamma, L} (1 + |\xi|)^{m - |\gamma|},\: \forall \alpha, \forall \gamma.$$
(ii) for $|\xi| \leq L$ we have\\
$$|\partial_x^{\alpha} \partial_{\xi}^{\gamma} a(x, \xi; h)| \leq C_{\alpha, \gamma, L} h^{-k - \delta(|\alpha| + |\gamma|)},\: \forall \alpha, \forall \gamma.$$

 For $a \in S_{\delta}^{k, m}(V),$ consider the operator
$$\Bigl({\rm Op}_h(a) f\Bigr)(x)= (2 \pi h)^{-d+1} \int\int e^{\ii \langle x - y, \xi \rangle/h} a(x, \xi; h) f(y) dy d\xi.$$
We have a calculus for the $h$-pseudo-differential operators with symbols in $S_{\delta}^{k, m}$ if $0 < \delta < 1/2.$
In particular, if $a \in S_{\delta}^{0, 1},\: b \in S_{\delta}^{0, -1},$ one gets
$$\|{\rm Op}_h(a) {\rm Op}_h(b) - {\rm Op}_h(a b) \|_{L^2} \leq C h^{1 - 2\delta}.$$
We refer to \cite{DS} and \cite{V1} for the calculus of $h$-pseudo-differential operators.\\

Let $D_{\nu} = - \ii \partial_{\nu}$, and let $\gamma_0$ denote the trace on $\Gamma.$ It is important to construct a semi-classical parametrix for the problem (\ref{eq:3.1}) in $Z_1 \cup Z_2 \cup Z_3$ and to find an approximation for the (exterior) {\it semi-classical} Dirichlet-to-Neumann map defined by 
\begin{equation} \label{eq:3.2}
{\mathcal N}_{ext}(z, h): H_h^s(\Gamma) \ni f \longrightarrow \gamma_0 hD_{\nu}u \in H_h^{s-1}(\Gamma).
\end{equation}
Here for $s \in \R,\:$ $H_h^s(\Gamma)$ is the semi-classical Sobolev space with norm $\|\langle h D\rangle^s u\|_{L^2(\Gamma)}$. Vodev \cite{V1} constructed a semi-classical parametrix $\tilde{u}$ where the equation in (\ref{eq:3.1}) is satisfied for $x \in K$. In fact the construction in \cite{V1} is made in a very small neighborhood of the boundary $\Gamma$ and the {\it local} parametrix is a Fourier integral operator with complex phase function. By using the resolvent $(-h^2 \Delta_D - z)^{-1}$ of the Dirichlet Laplacian in $\Omega$, one may modify the proof in \cite{V1} to obtain a parametrix in $\Omega$ (see \cite{P1} for more details).

To describe the local parametrix, consider {\it normal geodesic coordinates} $(x_1, x')$ in a neighborhood of a fixed point $x_0 \in \Gamma$, where $x_1 = {\rm dist} \:( x, \Gamma).$ Then locally the boundary $\Gamma$ is given by $x_1 = 0.$  Let $\psi(x') \in C_0^{\infty}(\Gamma)$ be a cut-off function with support in a small neighborhood of $x_0 \in \Gamma$ and $\psi(x') = 1$ in another neighborhood of $x_0.$ Then $-\frac{h^2}{c(x)}\nabla c(x) \nabla - z \frac{n(x)}{c(x)}$ in these coordinates has the form
$${\mathcal P}(z, h) = h^2D_{x_1}^2 + r(x, hD_{x'}) + q(x, h D_x) + h^2 \tilde{q}(x)- z m(x).$$
 with 
$$D_{x_1} = -\ii  \pa_{x_1},\: D_{x'} = -\ii \pa_{x'}, \: m(x) = \frac{n(x)}{c(x)},\: r(x, \xi') = \langle R(x)\xi', \xi'\rangle, \:q(x, \xi) = \langle q(x), \xi \rangle.$$
Here $R(x)$ is a symmetric $(d-1) \times (d-1)$ matrix with smooth real-valued enters and $r(0, x', \xi') = r_0(x', \xi')$ is the principal symbol of the Laplace-Beltrami operator $-\Delta_{\Gamma}$ on $\Gamma$.  Let 
$$\rho = \sqrt{z\:m(x) - r_0(x', \xi')} \in C^{\infty}(T^*(\Gamma))$$
 be the root of the equation
$\rho^2 + r_0(x', \xi') - z m(x) = 0$ with $\im \rho > 0$.
Let $\phi(\sigma) \in C^{\infty}(\R)$ be cut-off function such that $\phi(\sigma) = 1$ for $|\sigma| \leq 1,$ $\phi(\sigma) = 0$ for $|\sigma| \geq 2.$
In \cite{V1} for small $\delta_1 > 0$ and for $x$ close to the boundary it was constructed a parametrix 
\begin{equation} \label{eq:3.3}
\begin{cases}\tilde{u}_{\psi}(x) = ( 2 \pi h)^{-d +1} \int\int e^{\frac{\ii}{h} \varphi(x, y', \xi', z)}\phi\Bigl(\frac{x_1}{\delta_1}\Bigr)\\
\times\phi\Bigl(\frac{x_1}{\delta_1 \rho_1}\Bigr) a(x, \xi', z; h) f(y') dy'd\xi',\\
\tilde{u}_{\psi}\vert_{x_1 = 0} = \psi f,\end{cases}
\end{equation}
where $0 < \delta_1 < 1$ is small enough and $\rho_1 = 1$ if $ z \in Z_2 \cup Z_3,\: \rho_1 = |\rho|^3$ if $z \in Z_1$.  
 The phase $\varphi(x, y', \xi', z)$ is complex-valued and we have
$$\varphi\vert_{x_1 = 0}= - \langle x'- y', \xi'\rangle,\:\pa_{x_1} \varphi\vert_{x_1 =0} = \rho,\: \im \varphi \geq x_1 \im \rho/2,$$
while $a\big\vert_{x_1 = 0} = \psi(x').$ Next, $a = \sum_{k=0}^{N-1} \sum_{j=0}^{N-1} x_1^k h^j a_{k,j}(x', \xi', z),$
$$\varphi = - \langle x'- y', \xi'\rangle + \sum_{k=1}^{N-1} x_1^k \varphi_k(x', \xi', z),\: \varphi_1 = \rho,$$
$N \gg 1$ being a large integer. The phase $\varphi$ and the amplitude $a$ are determined so that
$$e^{-\frac{\ii\varphi}{h}} {\mathcal P}(z, h) e^{\frac{\ii\varphi}{h}}a = x_1^N A_N(x, \xi', z; h) + h^N B_N(x, \xi', z; h),$$
where $A_N,\: B_N$ are smooth functions and their behavior for $|\xi'| \to \infty$ is related to negative powers of $|\rho|.$ For example,
 $$|\pa_{x_1}^k \pa_{x'}^{\alpha} \pa_{\xi'}^{\beta} A_N (\phi(\delta_0 r_0(x', \xi'))| \leq C_{k, \alpha, \beta} |\rho|^{2-3N-3k - 2|\alpha|-2|\beta|}.$$
Moreover, for $x_1 > 0$ the parametrix $\tilde{u}_{\psi}$ has a decay ${\mathcal O}\Bigl(e^{- x_1\frac{|\im z| }{2 |\rho| h}}\Bigr)$ and for $ x_1 \geq |\rho|^3/\delta$ we get an estimate ${\mathcal O}\Bigl(e^{- C\frac{|\rho|^2|\im z| }{ h}}\Bigr)$.

Consider the (interior) semi-classical Dirichlet-to-Neumann map $\nc_{int}(z, h)f = \gamma_0 \partial_{\nu} u,$ related to the problem
\begin{equation} \label{eq:3.4} \begin{cases} (-\frac{h^2}{n(x)} \nabla c(x) \nabla - z )u =0 \:\: {\rm in}\: K,\\
u = f \:\: {\rm on}\: \Gamma, \end{cases}
\end{equation}
where $n(x) > 0, \: c(x) > 0$ are $C^{\infty}$ functions on $\Gamma$.
Then we have the following

\begin{proposition} [\cite{V1}] Given $0 < \epsilon \ll 1,$ there exists $0 < h_0(\epsilon) \ll 1$ such that for $z \in Z_1$ and $0 < h \leq h_0(\epsilon)$ we have
\begin{equation} \label{eq:3.5}
\|{\mathcal N}_{int}(z, h) f- {\rm Op}_h(\rho + h b)f\|_{H_s^1(\Gamma)} \leq \frac{Ch}{\sqrt{|\im z|}} \|f\|_{L^2(\Gamma)},
\end{equation} 
where $b \in S_0^{0, 0}(\Gamma)$ does not depend on $z, h$ and the function $n(x).$ Moreover, for $z \in Z_2\cup Z_3$ the above estimate holds with $|\im z|$ replaced by $1$.
\end{proposition}

With some modifications of the proof the same result remains true for unbounded domains $\R^d \setminus \bar{K}$ and obtain the estimate (\ref{eq:3.5}) for the semi-classical Dirichlet-to-Neumann operator ${\mathcal N}_{ext}(z, h)$ related to the problem (\ref{eq:3.1}) with $n(x) = c(x) = 1.$
(see \cite{P1}).

\section{Location of the eigenvalues of $G$}
\setcounter{equation}{0}
Let $u =(u_1, u_2) \not= 0$ be an eigenfunction of $G$ with eigenvalue $\lambda,\: \re \lambda< 0,$ and let $f = u_1\vert_{\Gamma}.$ Then from (\ref{eq:2.2}) we deduce
$(-\Delta + \lambda^2) u_1 = 0$ and
$\pa_{\nu} u_1 - \lambda \gamma u_1 = 0$ on $\Gamma$. Setting 
$$\lambda = \frac{\ii \sqrt{z}}{h},\: 0 < h \ll 1,$$
for $z \in Z_1 \cup Z_2 \cup Z_3,$
 one obtains the problem
$$\begin{cases} (- h^2 \Delta - z) u_1 = 0\: {\rm in }\: \:\Omega,\\
{\mathcal  N}_{ext}(z, h)f - \sqrt{z} \gamma f = 0\:\: {\rm on}\:\Gamma. \end{cases}$$

Consider the case (A) and notice that there exists $\ep_0 > 0$ such tat
$$0 < \epsilon_0 \leq \gamma(x) \leq 1 - \epsilon_0,\:\: \forall x \in \Gamma.$$ 
We will discuss the case $z \in Z_1$, the case $z \in Z_2 \cup  Z_3$ is more simple.
According to Proposition 1 for ${\mathcal N}_{ext}(z, h)$, for $z \in Z_1,\:1 \geq \im z \geq h^{\delta},\: \delta = 1/2 - \ep,$ we have
\begin{equation} \label{eq:4.1}
\|{\rm Op}_h(\rho)f  - \gamma \sqrt{z}f\|_{\lg} \leq C \frac{h}{\sqrt{|\im \:z|}} \|f\|_{\lg},
\end{equation}
while for $z \in Z_2 \cup Z_3$ the above estimate holds with $|\im\: z|$ replaced by 1.
Consider the symbol
$$c(x',\xi', z) = \rho(x',\xi', z)- \gamma\sqrt{z}=\frac{(1 - \gamma^2) z - r_0(x', \xi')}{\rho(x', \xi', z) + \gamma \sqrt{z}}. $$
We will show that $c(x', \xi', z)$ is elliptic in a suitable class.

Clearly, $c$ is elliptic for $|\xi'|$ large enough. So it remains to examine the behavior of $c$ for $|\xi'| \leq C_0$ and for these values of $\xi'$  we have $|\rho + \gamma \sqrt{z} | \leq C_1.$ 
 Introduce the set  
$${\mathcal F} = \{(x', \xi'):\:
 |1- r_0(x', \xi')| \leq \frac{\epsilon_0^2}{2}\}.$$
Then $\re \Bigl((1 - \gamma^2) z - r_0 \Bigr) = 1- r_0 - \gamma^2 \leq -\frac{\epsilon_0^2}{2}.$
If $(x', \xi') \notin {\mathcal F}$, we get
$$\im \Bigl( (1 - \gamma^2) z - r_0 \Bigr) = (1 - \gamma^2) \im z  \geq (1 - \gamma^2) h^\delta \geq \ep_1 h^\delta, \:\ep_1 > 0.$$
Consequently, the symbol $c$ is elliptic and
$$\im (\rho + \gamma \sqrt{z}) = \im \rho + \gamma \im \sqrt{z} \geq C h^{\delta}.$$ 
Thus, for bounded $|\xi'|$ we have $|c| \geq C_3 h^{\delta}, C_3> 0,$
while for large $|\xi'|$ we have $|c| \sim |\xi'|.$ Introduce the function
$$ \chi(x', \xi') = \phi(\delta_0 r_0(x', \xi')),\: 0 < \delta_0 \leq 1/2$$
and define ${\mathcal M}_1: = Z_1 \times {\rm supp}\, \chi,\: {\mathcal M}_2: = (Z_1 \times {\rm supp}\: (1 - \chi)) \cup ((Z_2\cup Z_3) \times T^* \Gamma).$
Set $\langle \xi'\rangle = ( 1 + |\xi'|)^{1/2}.$ It is easy to see that for $(z, x', \xi') \in {\mathcal M}_1,\: \im\: z \neq 0,$ we have
\begin{equation} \label{eq:4.2}
\big | \pa_{x'} ^{\alpha} \pa_{\xi'}^{\beta} \rho \big | \leq C_{\alpha, \beta} |\im\: z|^{1/2 - |\alpha| - |\beta|},\: |\alpha| + |\beta| \geq 1,
\end{equation}
$|\rho| \leq C,$ while for $(z, x', \xi') \in {\mathcal M}_2$ we have
\begin{equation} \label{eq:4.3}
\big | \pa_{x'} ^{\alpha} \pa_{\xi'}^{\beta} \rho \big | \leq C_{\alpha, \beta} \langle \xi' \rangle ^{1- |\beta|}.
\end{equation}
Thus, we conclude that $c =(\rho - \gamma \sqrt{z})\in S^{0, 1}_{\delta}.$ A similar analysis shows that $|\im\: z| c^{-1}\in S_{\delta}^{0, -1},$ while for $z \in Z_2 \cup Z_3$ we get $c^{-1} \in S_{\delta}^{0, -1}.$
Therefore
$$\|{\rm Op}_h(c^{-1}) g \|_{L^2(\Gamma)} \leq  C |\im\: z|^{-1} \|g\|_{L^2(\Gamma)}$$
and we deduce
$$\|{\rm Op}_h(c^{-1}){\rm  Op}_h(c) f \|_{L^2(\Gamma)} \leq C_1 \frac{h}{|\im\: z|^{3/2}} \|f\|_{L^2(\Gamma)}. $$
A more fine analysis (see \cite{P1}) shows that
$$\|{\rm Op}_h(c^{-1}) {\rm Op}_h(c) f - f \|_{L^2(\Gamma)} \leq C_2 \frac{h}{|\im z|^{2}} \|f\|_{L^2(\Gamma)}.$$
Consequently,  one concludes that
\begin{equation} \label{eq:4.4}
\|f\|_{L^2(\Gamma)} \leq C_3 \Bigl(  h^{1-2\delta} + h^{1-\frac{3}{2}\delta}\Bigr) \|f\|_{L^2(\Gamma)}.
\end{equation}
Since $\delta = 1/2 - \epsilon,\: 0 < \epsilon \ll 1$, from (\ref{eq:4.4}) it follows $f = 0$ for $0 < h \leq h_0(\epsilon)$ small enough. Since $-h^2 \Delta$ with Dirichlet boundary conditions does not have eigenvalues in $\{z\in \C:\: \re \: z < 0\}$, one gets $u_1 = 0.$
Going back to the eigenvalues and using the scaling, one obtains that in the case (A) the eigenvalues of $G$ lie in the region
$$\Lambda_{\epsilon} = \{ z \in \C:\: |\re\: z | \leq C_{\epsilon} (|\im\: z|^{\frac{1}{2} + \epsilon} + 1), \: \re\: z < 0\}.$$

In the case (B) the above  analysis  works only for $z \in Z_1 \cup Z_3$.
Indeed for $z \in Z_1$  we have 
$$ \re ((1 - \gamma^2)- r_0) \leq (1 -\gamma^2) < - \eta_0 < 0$$ 
and  again $c \in S_{\delta}^{0, 1}, c^{-1} \in S_{\delta}^{0, -1}.$  Thus for $z \in Z_1 \cup Z_3$ we obtain that the eigenvalues $\lambda = \frac{\ii \sqrt{z}}{h}$ must lie in $\Lambda_{\epsilon}.$
For $z \in Z_2$ the argument, exploited in the case (A), breaks down since for $\re\: z = -1, \:\im\: z = 0$ the symbol 
$$ 1 + r_0(x', \xi') - \gamma(x') $$
is not elliptic and it may vanish for some $(x_0', \xi_0').$

 Let $z = - 1 + \ii \im\: z \in Z_2$. For such $z$ we have an better approximation $T(z, h)$ of the operator ${\mathcal N}_{ext}(z, h)$ (see \cite{V1}, \cite{P1}) for which we have
\begin{equation} \label{eq:4.5}
\|{\mathcal N}_{ext}(z, h)f - T(z, h) f\|_{H^{1}(\Gamma)} \leq C_N h^{-s_d + N} \|f\|_{\lg},\: \forall N \in \N,
\end{equation}
with $s_d > 0$ depending only on  the dimension $d$.
 Therefore, if $f$ is related to the trace of an eigenfunction of $G$, from the equality ${\mathcal N}_{ext}(z, h)f - \gamma\sqrt{z}f = 0$ on $\Gamma$ we obtain
$$|\re\: ( T(z, h) f - \gamma \sqrt{z} f, f)_{L^2(\Gamma)} | \leq C_N h^{-s_d + N} \|f\|_{\lg}.$$
Next, by applying Taylor formula, we write
\begin{eqnarray} \label{eq:4.6}
\re \:\Bigl( (T(z, h) - \gamma \sqrt{z})f, f\Bigr)_{\lg} = \re \Bigl((T(-1, h) - \ii \gamma)f, f\Bigr)_{\lg}\nonumber\\
 - \im\: z \im \Bigl(\Bigl[ \frac{\partial T}{\partial z} (z_t, h) - \gamma \frac{1}{2 \sqrt{z_t}}\Bigr]f, f\Bigr)_{\lg}
\end{eqnarray}
with $z_t = - 1 + \ii\: t \im\: z \in Z_2, \: 0 < t < 1.$ We may replace in (\ref{eq:4.6}) the operator $\frac{\partial T}{\partial z}(z_t, h)$ by the operator
${\rm Op}_h(\frac{d\rho}{dz}(z_t, h))$
 modulo $\oc(h)\|f\|_{\lg}^2 $ term and a sharp analysis shows that 
$$\im \Bigl( \Bigl({\rm Op}_h(\frac{d\rho}{dz}(z_t, h)) - \gamma \frac{1}{2 \sqrt{z_t}}\Bigr) f, f \Bigr)_{L^2(\lg)} \geq \alpha_0 \|f\|_{\lg}^2, \: \alpha_0 > 0.$$
We refer to \cite{P1} for the details of this argument. Combining (\ref{eq:4.5}) and (\ref{eq:4.6}), one estimates $|\im z|$ and for small $h$ and every $N \in \N$, we obtain that the eigenvalues $\lambda = \frac{\ii \sqrt{z}}{h}$ of $G$ with $z \in Z_2$ must lie in the region
$${\mathcal R}_N = \{z \in \C:\: |\im\: z| \leq C_N (|\re\: z| + 1)^{-N},\: \re\: z < 0\}.$$
Finally, we have the following
\begin{theorem}[\cite{P1}] In the case $(A)$ for every $\epsilon, \: 0 < \epsilon \ll 1,$ the eigenvalues of $G$ lie in the region $\Lambda_{\epsilon}.$ In the case $(B)$ for every  $\epsilon, \: 0 < \epsilon \ll 1,$ and every $N \in \N$ the eigenvalues of $G$ lie in the region $\Lambda_{\epsilon} \cup {\mathcal R}_N$.
\end{theorem}

For {\it strictly convex} obstacles $K$ we have a more precise result concerning the operator ${\mathcal N}_{out}(z, h)$ based on the construction of a semi-classical parametrix for the problem (\ref{eq:3.1}) when $\re z = 1$ and $h^{1/2 - \epsilon} \geq \im z \geq h^{2/3}$ (see  \cite{V2}, \cite{P1}) or  $0 < \im z \leq h^{2/3}$  (see \cite{Sj}). This makes possible to improve the above result in the case (B) and to obtain the following
\begin{theorem}[\cite{P1}]In the case $(B)$ for  every $N \in \N$  outside  the region $ {\mathcal R}_N$ we have only finite number eigenvalues of the generator $G$. 
\end{theorem}

Moreover, we have the following

\begin{proposition}[\cite{CPR1}] Assume that $d$ is odd. Then the operator $G$ has no a sequence of eigenvalues $\lambda_j,\, \re \lambda _j < 0$ such that $\lim_{j \to \infty} \lambda_j = \ii z_0,\, z_0 \in \R.$
\end{proposition}

It is world noting that the Dirichlet-to-Neumann map can be used to establish the discreetness of the spectrum of $G$ in $\{z \in \C: \re\: z < 0\}$. We follow below the argument of \cite{P1}.  For  $\re\: \lambda < 0$ introduce the map
$$\nc(\lambda): H^s(\Gamma) \ni f \longrightarrow \pa_{\nu} u\vert_{\Gamma} \in H^{s-1}(\Gamma),$$
where $u$ is the solution of the problem 
\begin{equation} \label{eq:4.7}
\begin{cases} (\Delta - \lambda^2) u = 0 \:\: {\rm in}\:\: \Omega,\: u \in H^2(\Omega),\\
u = f \:{\rm on}\: \Gamma. \end{cases}
\end{equation}
It is well known that $\nc(\lambda)$ is a meromorphic function in $\C$ for $d$ odd and in the logarithmic covering of $\C$ for $d$ even and the poles of $\nc(\lambda)$ in $\C \setminus \{0\}$ coincide with the resonances of the Dirichlet problem for the Laplacian (see for instance, \cite{Sj}). On the other hand,  $u \in H^2(\Omega)$ implies that $u$ is $\lambda-$incoming in the sense of Lax and Phillips (see Chapter IV in \cite{LP1}). 
Notice that the definition of outgoing/incoming solutions in \cite{Sj} is different from that in \cite{LP1} and the resonances in \cite{Sj} lie in $\im\: z < 0$, while in \cite{LP1} they are in the half-plan $\im\: z > 0.$ Consequently, $\nc(\lambda)$ is analytic for $\re\: \lambda < 0.$  The same is true for the Neumann problem for the $\Delta - \lambda^2$, hence $\nc^{-1}(\lambda)$ is also analytic for $\re\: \lambda < 0$ and the poles of $\nc^{-1}(\lambda)$ are the resonances of the Neumann problem (\cite{SjV}). Therefore, the boundary condition in (\ref{eq:2.2}) may be written as follows
$$\nc(\lambda) \Bigl( I - \lambda \nc^{-1}(\lambda) \gamma \Bigr) f_1 = 0,\:\re \:\lambda < 0,\: x \in \Gamma.$$
The operator $\nc(\lambda): L^2(\Gamma) \longrightarrow H^{1}(\Gamma)$ is compact and Theorem 1 guarantees that there are points $\lambda_0,\: \re\: \lambda_0 < 0,$ for which  $(I - \lambda_0 \nc^{-1}(\lambda_0)\gamma)$ is invertible. Applying the analytic Fredholm theorem, we  conclude that the spectrum of $G$ in $\{z \in \C:  \re \:\lambda < 0\}$ is formed by isolated eigenvalues with finite multiplicities.\\

We finish this section by a trace formula involving the operator 
$$ C(\lambda) := \nc(\lambda) - \lambda \gamma = \nc(\lambda) \Bigl( I - \lambda \nc^{-1}(\lambda) \gamma \Bigr),$$
which is an analytic operator-valued function in $\{ z \in \C: \: \re \: z < 0\}$, while $C(\lambda)^{-1}$ is meromorphic in the same domain.
 We wish to find  a formula for the trace
\begin{equation} \label{eq:4.8}
  {\rm tr}\: \frac{1}{2 \pi i} \int_{\delta} (\lambda - G)^{-1} d\lambda,
\end{equation}
 where $\omega \subset \{\re\: z < 0\}$ has as a boundary  the curve $\delta$ and $(G-\lambda)^{-1}$ is analytic on $
\delta$. We know that $(G- \lambda)^{-1}$ is meromorphic in $\omega$ and if $\lambda_0$ is a pole of $(G- \lambda)^{-1},$ then the multiplicity of the eigenvalue $\lambda_0$ of $G$ is given by
$${\rm mult}\: (\lambda_0) ={\rm rank} \frac{1}{2 \pi i} \int_{|\lambda - \lambda_0| = \ep_0} (\lambda - G)^{-1} d\lambda,$$
with $\ep_0 > 0$ small enough. Therefore, (\ref{eq:4.8}) is equal to the number of the eigenvalues of $G$ in $\omega$ contented with their multiplicities.

Let $(u, w) =  (G-\lambda)^{-1} (f, g)$. Then $w = \lambda u + f$ and setting $q = u\vert_{\Gamma}$, one gets
$$u = R_D(\lambda) (g + \lambda f) + K(\lambda) q.$$
Here $R_D(\lambda) = (\Delta_D - \lambda^2)^{-1} $ is the resolvent of the operator $\Delta_D$ with Dirichlet boundary conditions
 and $K(\lambda)$ satisfies
$$\begin{cases}  (\Delta - \lambda^2) K(\lambda) = 0 \:\: {\rm in}\: \Omega,\\
K(\lambda) = \:Id \:\:{\rm on}\: \Gamma.
\end{cases}$$
The boundary condition on $\Gamma$ yields
$$\pa_{\nu} [R_D(\lambda) ( g + \lambda f) + K(\lambda) q]- \gamma \lambda [ R_D(\lambda) ( g + \lambda f) + K(\lambda) q] - \gamma f = 0,\: x \in \Gamma$$
and the term $\gamma \lambda [R_D(\lambda) ( g + \lambda f) $ vanishes. Since $\nc(\lambda) = \pa_{\nu} K(\lambda)\vert_{\Gamma}$ is the Dirichlet-to-Neumann map, assuming that $C^{-1}(\lambda)$ is invertible, 
we deduce
 $$q = C^{-1}(\lambda) \Bigl( [ \pa_{\nu} R_D(\lambda)(g + \lambda f) ] - \gamma f \Bigr).$$
Therefore
$$u = \Bigl[\lambda R_D(\lambda) + K(\lambda) C^{-1}(\lambda) \lambda \pa_{\nu} R_{D}(\lambda) - \gamma\Bigr] f + X g,$$
$$ w = Y f + \Bigl[\lambda R_{D}(\lambda) + \lambda K(\lambda) C^{-1}(\lambda) \pa_{\nu} R_D(\lambda) \Bigr]g,$$
where the operators $X$ and $Y$ are not important for the calculus of the trace.
Thus we are going to study the integral
$${\rm tr}\: \int_{\delta} \Bigl(2 \lambda K(\lambda) C^{-1}(\lambda) \pa_{\nu} R_D(\lambda) - C^{-1}(\lambda) \gamma\Bigr) d\lambda.$$
For the first term  we apply the cyclicity of the trace and the fact that
$$\frac{\pa \nc}{\pa \lambda}(\lambda)= \pa_{\nu} \frac{\pa K}{\pa \lambda}(\lambda)  = 2 \lambda \pa_{\nu}R_{D} (\lambda) K(\lambda).$$
Finally, we obtain the following
\begin{proposition}[\cite{P1}] Let $\delta \subset \{ z \in \C: \: \re\: \lambda < 0\}$ be a closed positively oriented curve and let $\omega$ be the domain bounded by $\delta$. Assume that  $C^{-1}(\lambda)$ is meromorphic in $\omega$ without poles on $\delta$ . Then
\begin{equation} \label{eq:4.9}
 {\rm tr}\: \frac{1}{2 \pi i} \int_{\delta} (\lambda - G)^{-1} d\lambda = {\rm tr} \frac{1}{2 \pi i} \int_{\delta} C^{-1}(\lambda) \frac{\pa C}{\pa \lambda}(\lambda) d \lambda.
\end{equation}
\end{proposition}

The idea to write the right-hand side of (\ref{eq:4.9}) as the trace of an integral involving the product of a meromorphic function $T^{-1}(\lambda)$ and its derivative $\frac{d T}{d \lambda}(\lambda)$ going back to \cite{SjV}, \cite{CPV} (see also Proposition 3 in the next section).
We expect that in the case (B) Proposition 3 combined with the techniques in \cite{SjV} will imply a Weyl formula for the eigenvalues of $G$ lying in ${\mathcal R}_N.$

We conjecture that for $N$ large enough and $\gamma(x) > 1,\, \forall x \in \Gamma,$ the counting function 
$$N(r) = \#\{ \lambda_j \in \sigma_p(G):\, |\lambda_j| \leq r, \,\, \lambda_j\in {\mathcal R}_N\} $$ 
has the asymptotic
\begin{eqnarray} \label{eq:4.10}
N(r) = (2\pi)^{-d+1} \omega_{d-1} \Bigl(\int_{\Gamma} (\gamma^2(y') -1)^{(d-1)/2} dy' \Bigr)r^{d-1}\nonumber \\
 + {\mathcal O}_{\gamma}(r^{d-2}),\: r \geq r_0(\gamma),
\end{eqnarray}
where $\omega_{d-1} = {\rm vol}\: \{x \in \R^{d-1}: |x| \leq 1\}.$ For strictly convex obstacles and $\gamma(x) > 1$ this will imply a Weyl asymptotics of all eigenvalues of $G$. Notice that for ball $B_3$ we have the following 

\begin{proposition} [\cite{P1}] For $\gamma \equiv const > 1$ and $K = B_3$ all eigenvalues $\lambda_j$ of $G$ are real and they lie in the interval $(-\infty, -\frac{1}{\gamma - 1}].$
  Moreover, there is an infinite number  of  real eigenvalues of $G$.
\end{proposition}
Hence in this case we must study the asymptotic of $N(r)$ for $r \geq -\frac{1}{\gamma - 1} = r_0(\gamma).$ Moreover, following the
analysis in \cite{P1}, we may prove that (\ref{eq:4.10}) holds for $K = B_3$ and constant $\gamma$.

By a similar argument we may study the eigenvalues of the generator $G$ of the contraction semigroup associated to Maxwell's equations with dissipative boundary conditions

\begin{equation}  \label{eq:4.11}
\begin{aligned} 
&\partial_t E = \curl B,\:\:\partial_t B = -\curl E \: {\rm in}\:\:\R_t^+ \times \Omega,
\\
&E_{tan} - \gamma(x)(\nu(x) \wedge B_{tan}) = 0 \:\:{\rm on} \:\R_t^+ \times \Gamma,
\\
&E(0, x) = e_0(x), \:\: B(0, x) = b_0(x). 
\end{aligned}
\end{equation}

The solution of the problem (\ref{eq:4.1}1) is given by a contraction semigroup 
$$(E, B) = V(t)f = e^{tG_b} f,\: t \geq 0,$$
 where the generator
$G_b$ has  domain $D(G_b)$
that is the closure in the graph norm of functions $u = (v, w) \in (C_{(0)}^{\infty} (\R^3))^3 \times (C_{(0)}^{\infty} (\R^3))^3$ satisfying the boundary condition $v_{tan} - \gamma (\nu \wedge w_{tan}) = 0$ on $\Gamma.$ Here $u_{tan} = u - \la u, \nu \ra \nu.$ 
For Maxwell's equations for $0 < \gamma(x) < 1$ and $\gamma(x) > 1$ we have the same location of eigenvalues of $G_b$. This location has been examined  in \cite{CPR2} by a semi-classical analysis of a $h$-pseudo-differential system on the boundary $\Gamma$. We have the following
\begin{theorem}[\cite{CPR2}] Assume that for all $x\in \Gamma$, $\gamma(x) \neq 1$.
Then for every $0 < \ep \ll 1$ and every $N \in \N$ there 
are constants $C_{\ep} > 0$ and $C_N > 0$ such that the eigenvalues of $G_b$ lie in the region $\Lambda_{\ep} \cup {\mathcal R}_N$, where
$$\Lambda_{\ep} = \{ z \in \C: \: |\re z | \leq C_{\ep} (|\im z|^{1/2 + \epsilon} + 1),\: \re z < 0\},$$
$${\mathcal R}_N = \{ z \in \C: \: |\im z| \leq C_N ( |\re z| + 1)^{-N},\: \re z < 0\}.$$
\end{theorem}

It is interesting to notice that for Maxwell's equation if $\gamma(x) \equiv 1,\, \forall x \in\Gamma,$ and $K = B_3$ is the unit ball in $\R^3$, then $G_b$ has no eigenvalues (see \cite{CPR2} for other results concerning the case $\gamma = const$ and $B_3$).

\section{Location and Weyl formula for the (ITE)}
\renewcommand{\theequation}{\arabic{section}.\arabic{equation}}
\setcounter{equation}{0}

To examine the location of the (ITE), set $\lambda = \frac{z}{h^2},\: z \in Z_1 \cup Z_2 \cup Z_3$. If $\lambda$ is an (ITE) with eigenfunction $(u, w)$, consider $u\vert_{\Gamma} = w\vert_{\Gamma} = f.$ Introduce the Dirichlet-to-Neumann operators $\nc_j = \nc_j(z, h),\: j = 1, 2$ related to 
$${\mathcal P}_j(z, h) = -\frac{h^2}{n_j(x)} \nabla c_j(x) \nabla - z \frac{c_j(x)}{n_j(x)},\: j = 1,2.$$
The boundary condition in the problem (\ref{eq:2.5}) implies
$$c_1 \nc_1(z, h) f - c_2 \nc_2(z, h) f = 0.$$
As in the Section 3, one introduces normal geodesic coordinates $(x_1, x')$ and define
 $$\rho_j = \sqrt{z\: \frac{n_j(x)}{c_j(x')} - r_0(x', \xi') }, j = 1, 2$$
with $\im \: \rho_j > 0.$ Applying Proposition 1 for the operators $\nc_j(x, h)$,  we deduce
$$\|c_1 {\rm Op}_h(\rho_1)f - c_2 {\rm Op}_h(\rho_2)f\|_{\lg} \leq \frac{C h}{\sqrt{|\im\: z|}}\|f\|_{\lg}.$$
Below we discuss only the case $c_1(x) = c_2(x) \equiv 1,\,\forall x \in \Gamma$. Then we have a better estimate
\begin{equation} \label{eq:5.1}
\|{\rm Op}_h(\rho_1)f -  {\rm Op}_h(\rho_2)f\|_{H^1_h(\Gamma)} \leq \frac{C h}{\sqrt{|\im\: z|}}\|f\|_{\lg}
\end{equation}
and we must invert the operator ${\rm Op}_h(\rho_1) -  {\rm Op}_h(\rho_2)$. Writing
$$\rho_1 - \rho_1 =   \frac{z(n_1(x') - n_2(x'))}{\rho_1 + \rho_2},$$
it is easy to see that $\rho_1 - \rho_2$ is elliptic and $(\rho_1 - \rho_2)^{-1} \in S_{\delta}^{0, -1}$ for $z \in Z_1$, while $(\rho_1 - \rho_2)^{-1} \in S^{0, -1}$ for $z \in Z_2 \cup Z_3.$ For $\delta = 1/2-\ep < 1/2$ we may use the calculus of h-pseudo-differential operators and (\ref{eq:5.1}) implies, as in Section 4,  $f = 0.$ The latter yields $u = w = 0.$ Returning to the eigenvalues $\lambda = \frac{z}{h^2}$, we get that the (ITE) lie in the domain $\Lambda_{+}$ defined below.
The analysis of the general case when $c_j(x)$ are not equal to 1 is more complicated and we refer to \cite{V1} for the details. Thus we have the following

\begin{theorem}[\cite{V1}] Assume $(\ref{eq:2.6})$ fulfilled together with the condition 
$$c_1(x) = c_2(x), \: \partial_{\nu} c_1(x) = \partial_{\nu} c_2(x), \: \forall x \in \Gamma.$$
 Then for every $0 < \epsilon \ll 1$ the (ITE) lie   the region
$$\Lambda_{+, \ep}: = \{z \in \C:\: \re\: \lambda \geq 0, \: |\im \:\lambda| \leq C_{\epsilon}(\re\: \lambda + 1)^{3/4 + \epsilon}\}$$ 
and there are only a finite number $(${\rm ITE}$)$ with $\re\: \lambda < 0.$ If $c_1(x) \neq c_2(x), \forall x \in \Gamma,$ the $(${\rm ITE}$)$\: lie in
$$\Lambda_{+, \ep}': = \{z \in \C:\: \re \:\lambda \geq 0, \: |\im\: \lambda| \leq C_{\epsilon}(\re\: \lambda + 1)^{4/5 + \epsilon}\}.$$ If $(c_1(x) - c_2(x))d(x) > 0,\: \forall x \in \Gamma,$ we have only a finite number $(${\rm ITE}$)$ with $\re\:\lambda < 0$. 
Moreover, if we assume that $(c_1(x) - c_2(x))d(x) < 0,\: \forall x \in \Gamma,$ then for $\re\: \lambda \geq 0$ the $(${\rm ITE}$)$ are in  $\Lambda_{+}$, while for $\re\: \lambda < 0$ and every $N \geq 1$ there exists $C_N > 0$ such that $(${\rm ITE}$)$\: lie in
$${\mathcal R}_N = \{\lambda \in \C:\: |\im\: \lambda| \leq C_N (|\re\: \lambda| + 1)^{-N},\: \re\: \lambda \leq 0\}.$$
\end{theorem}
A weaker result in a partial case $n_1(x) \equiv 1, n_2(x) > 1$ in $K$ with eigenvalues-free region
$$ \{z \in \C:\: \re\: \lambda \geq 0, \: |\im\: \lambda| \geq C(\re \:\lambda + 1)^{24/25}\}$$ 
has been obtained in \cite{HK}.

For strictly convex obstacles one may construct a parametrix for the problem (3.4) and $\re z = 1,\,h^{1/2 - \epsilon} \geq \im z \geq h^{1- \epsilon}$ by using more complicated construction and exploiting the properties of the Airy function Ai(z) (see \cite{V2} for more details). 
 This leads to the following improvement of Theorem 4.

\begin{theorem} [\cite{V2}] Assume $K$ strictly convex, the condition $(\ref{eq:2.6})$ satisfied and $c_1(x) = c_2(x), \pa_{\nu} c_1(x) = \pa_{\nu}c_2(x), \: x\in \Gamma.$ Then for every $\ep > 0$
the $(${\rm ITE}$)$ lie in the region
$$M_{+, \ep}: = \{z \in \C:\: \re\: \lambda \geq 0, \: |\im\: \lambda| \leq C_{\epsilon}(\re \:\lambda + 1)^{1/2 + \epsilon}\}$$ 
and there are only a finite number $(${\rm ITE}$)$ with $\re \:\lambda < 0.$
\end{theorem}

This results is almost optimal, since for the unit ball in $\R^d$ we have the following

 \begin{theorem} [\cite{PV1}]Let $K =\{x\in \R^d:\,|x|\le 1\},\: d \geq 2$. Suppose that the functions $c_j$, $n_j$, $j=1,2$, are constants everywhere in $K$, $c_1 = c_2,$ and the condition $(\ref{eq:2.6})$ is satisfied. 
Then, there are no $({\rm ITE})$ in the region ${\mathcal M}_{+, 0}$
\end{theorem}
The case $d =1$ and $K = \{ x \in \R:\: |x| \leq 1\}$ has been previously examined in \cite{S} and \cite{PS}.\\

Now we pass to the Weyl formula for the counting function $N(r)$ of the (ITE) and introduce the coefficients
$$\tau_j = \frac{\omega_d }{(2\pi)^d} \int_K \Bigl(\frac{n_j(x)}{c_j(x)}\Bigr)^{d/2} dx,\: j = 1, 2,$$
where $\omega_d$ is the volume of the unit ball in $\R^d.$

In the anisotropic case $c_1(x)= 1, n_1(x) = 1, c_2(x) \neq 1, c_2(x) n_2(x) \neq 1, \forall x \in \bar{K},$ the asymptotics
\begin{equation} \label{eq:5.2}
N(r) \sim (\tau_1 + \tau_2) r^d,\: r \to +\infty.
\end{equation}
has been obtained by Lakshatanov and Vainberg \cite{LV} under some additional
assumptions which guarantee that the boundary problem is {\it parameter-elliptic.} 

By the results of Agranovich and Vishik \cite{AV} for the closed operator ${\mathcal A}$ related to (\ref{eq:2.5}) outside every angle $D_{\alpha} = \{z \in \C:\: |\arg z| \leq \alpha\}$, we have only a finite number of (ITE) and  the following estimate holds
$$\|(z - {\mathcal A})^{-1}\| \leq C_{\alpha} |z|^{-1},\: z \notin D_{\alpha}, \: |z| \gg 1.$$
The authors applied directly a result of Boimanov-Kostjuchenko \cite{BK} leading to (\ref{eq:5.2}).\\

The isotropic case $c_1(x)  = c_2(x) = 1,\: \forall x \in \bar{K},$  $n_1(x) = 1,\: n_2(x) \neq 1, \: \forall x \in \Gamma,$ is more difficult since the corresponding operator ${\mathcal A}$ has domain
$$D({\mathcal A}) = \{ (u, w) \in L^2(K) \times L^2(K):\: \Delta u \in L^2(K),\: \Delta v \in L^2(K),$$
$$ u- w = 0, \partial_{\nu} (u - w) = 0\: {\rm on}\:\: \Gamma\}. $$
Thus  $D(A)$ is not included in $ H^2(K)$, and the problem is not parameter-elliptic. In this case Robbiano \cite{R} obtained (\ref{eq:5.2}) by establishing the asymptotics
$$\sum_{j} \frac{1}{|\lambda_j|^p + t} = \alpha t^{-1 + \frac{d}{2p}} + o(t^{-1 + \frac{d}{2p}}), \:t \to + \infty.$$
where $p \in \N$ is sufficiently large. An application of the  Tauberian theorem of Hardy-Littlewood yields the result. By this argument one obtains a very week estimate for the remainder which can be estimated by the principal term divided by a logarithmic factor. To get better results, it is important to take into account {\it parabolic eigenvalues-free regions} and to apply different techniques which are not based on Tauberian theorems.

\begin{theorem} [\cite{PV}]
Under the condition $(\ref{eq:2.6})$, assume that there are no $(${\rm ITE}$)$ in the region
\begin{equation} \label{eq:5.3}
\{\lambda \in \C:\: |\im \lambda| \geq C (|\re \lambda| + 1)^{1 - \frac{\kappa}{2}}\},\: C > 0, 0 < \kappa \leq 1.
\end{equation}
Then for every $0 < \epsilon \ll 1$ we have the asymptotics
\begin{equation} \label{eq:5.4}
N(r) = (\tau_1 + \tau_2) r^d + {\mathcal O}_{\epsilon} ( r^{d - \kappa + \epsilon}),\: r \to +\infty.
\end{equation}
\end{theorem}
$\:\:\:\:\:\:\bullet$ According to Theorem 4, for arbitrary obstacles and $c_1(x) = c_2(x),\: \partial_{\nu} c_1(x) = \partial_{\nu} c_2(x), \forall x \in \Gamma,$  we can take $\kappa = \frac{1}{2} - \epsilon$  and we obtain a remainder ${\mathcal O}_{\epsilon}(r^{d- 1/2 + \epsilon}).$

$\bullet$ Taking into account Theorem 5, for strictly convex obstacles we choose $\kappa = 1- \epsilon,\: \forall \epsilon$. Consequently, we have in this case a remainder $O_{\epsilon} (r^{d- 1 + \epsilon}).$

$\bullet$ The optimal result should be to have a eigenvalues-free region with
$\kappa = 1$ as it was proved in \cite{PV1}, \cite{S}, \cite{PS} for the case when $K$ is a ball and the functions $c_j, n_j$ are constants. However, even with $\kappa = 1$, to obtain an optimal remainder ${\mathcal O}(r^{d-1})$ some extra work is needed and this is an interesting open problem.

The proof of Theorem 7 is long and technical. After a semi-classical scaling, the idea is to reduce the analysis of $N(r)$ to the trace of an integral involving the product of a meromorphic function $T^{-1}(\lambda)$ and its derivative $\frac{dT}{d\lambda}(\lambda)$ similar to Proposition 3. 
 Set
$Z =  \{z \in \C;\:\frac{1}{2} \leq  |\re\: z| \leq 3, \: |\im\: z| \leq 1\}$
and consider for $z \in Z$ and $0 < h \ll 1$ the operator
$$h T(z/ h^2): = c_1 {\mathcal N}_1(z, h) - c_2{\mathcal N}_2(z, h),$$
where the DN-maps $N_j(z, h)$ are defined in the beginning of this section.

 Let  $G_D^{(j)},\: j = 1, 2,$ be the Dirichlet self-adjoint realization of the operator $L_j : = -n_j^{-1} \nabla c_j \nabla$ in the space $H_j = L^2( K, n_j(x) dx).$ Set ${\mathcal H} = H_1 \oplus H_2$ and let $R(\lambda)$ be the resolvent of the transmission boundary problem. 
We omit in the notation $j = 1,2$ and consider the operators
$${\mathcal N}(z, h) Op_h(1 -\chi)f = \tilde{{\mathcal N}}(z, h) f - \gamma_0 D_{\nu} (h^2 G_D - z) ^{-1}\frac{c}{n} {\rm Op}_h(p)f,$$
$$F(z, h) = {\mathcal N}(z, h) - \tilde{{\mathcal N}}(z, h) = {\mathcal N}(z, h) {\rm Op}_h(\chi) - \gamma_0 D_{\nu} (h^2 G_D - z) ^{-1} \frac{c}{n} {\rm Op}_h(p),$$
where $\chi(x', \xi') = \Phi (\delta_0 r_0(x', \xi'))$ with $\Phi(\sigma) = 1$ for $|\sigma| \leq 1$ and $\Phi(\sigma) = 0$ for $|\sigma| \geq 2$, while $0 < \delta_0 \ll 1$ is small enough. Here ${\tilde{\mathcal N}}(z, h)$ is the parametrix of the DN operator ${\mathcal N}(z, h){\rm Op}_h(1 - \chi)$ in the domain where $r_0(x', \xi') > \frac{1}{\delta_0}$ and $p$ is some symbol having  behavior ${\mathcal O}(h^N)$ with all its derivatives. The number $N$ will be taken large enough and it depends only on  the parametrix construction.\\

 The operator $F(z, h)$ is meromorphic with values in the space of trace class operators and we denote by $\mu_j(F(z, h))$ its characteristic eigenvalues.

\begin{lemma} If $z/h^2$ does not belong to spec $G_D$, then for every integer $0 \leq m \leq N/4$ we have
$$ \mu_j(F(z, h)) \leq \frac{C}{\delta(z, h)} \Bigl(h j ^{1/(d-1)}\Bigr)^{-2m},\: \forall j \in \N,$$
where $\delta(z, h): = \min\{1, {\rm dist}\: \{z, {\rm spec}\: h^2 G_D\}\} > 0$ and $C > 0$ depends on $m$ and $N$ but is independent of $z, h, j$.
\end{lemma}
Let
$$T(\lambda): = \gamma_0 c_1  D_{\nu} K_1(\lambda) - \gamma_0 c_2 D_{\nu} K_2(\lambda),$$
 where $K_j(\lambda) f = u,$ and $u$ is the solution of the problem
$$\begin{cases}   \Bigl( L_j - \lambda\Bigr) u = 0 \:\: {\rm in}\: K,\\
u = f \: {\rm on}\: \Gamma. \end{cases}.$$

\begin{proposition} Assume that $T(\lambda)^{-1}$ is a meromorphic function with residues of finite rank. Let $\delta \subset \C$ be a simple closed positively oriented curve which avoids the eigenvalues of 
$G_D^{(j)}$, $j=1,2$, as well as the poles of $T(\lambda)^{-1}$. Then we have the identity
\begin{eqnarray} \label{eq:5.5} -{\rm tr}_{{\cal H}}\, (2\pi i)^{-1}\int_{\delta} R(\lambda)d\lambda \nonumber =\sum_{j=1}^2{\rm tr}_{H_j}\, (2\pi i)^{-1}\int_{\delta}
(\lambda - G_D^{(j)})^{-1}d\lambda\nonumber \\
-{\rm tr}_{L^2(\Gamma)}\, (2\pi i)^{-1}\int_{\delta}T(\lambda)^{-1}\frac{dT(\lambda)}{d\lambda}d\lambda.
\end{eqnarray}
\end{proposition}

Let us mention that if $R(\lambda)$ is an operator-valued meromorphic function with residues of finite rank, the multiplicity of a pole
$\lambda_k \in \C$ of $R(\lambda)$ is defined by
$${\rm mult}\: (\lambda_k) = -{\rm rank}\: (2 \pi i)^{-1} \int_{|\lambda - \lambda_k| = \epsilon} R(\lambda) d\lambda,\: 0 < \epsilon \ll 1.$$
On the other hand, the rank of the operator above is equal to the trace of this operator and on the left-hand side of (\ref{eq:5.5}) we have the sum of the mutiplicities
of the (ITE) lying in the domain $\omega_{\delta}\subset  \C$ bounded by $\delta$. Clearly, the terms with $(\lambda - G_D^{(j)})^{-1}$ yield the sum of eigenvalues of $G_D^{(j)}$ in $\omega_{\delta}$ counted with their multiplicities.\\

It is possible to construct invertible, bounded operator
$E(z, h): H^s_h(\Gamma) \to H_h^{s+ 1}(\Gamma)$ with bounded inverse $E(z, h)^{-1}: H_h^s(\Gamma) \to H_h^{s-1}(\Gamma),\: \forall s \in \R,$ so that 
$$h T(z/ h^2) = E^{-1}(z, h) (I + {\mathcal K}(z, h)),$$
$$ (h T(z/ h^2))^{-1} = (I + {\mathcal K}(z, h))^{-1} E(z, h)$$
with a trace class operator 
$${\mathcal K}(z, h) = E(z, h) (c_1 F_1(z, h) - c_2 F_2(z, h)) + {\mathcal L}(z, h).$$
Moreover, the operators $E(z, h), E^{-1}(z, h),$ are holomorphic with respect to $z$ in  $Z$
while
${\mathcal K}(z, h)$ is meromorphic operator-valued function in this region.
Then
$${\rm tr}\: \int_{\delta} T^{-1}(z/h^2) \frac{d}{dz}T(z/h^2)dz = {\rm tr}\: \int_{\delta} (I + {\mathcal K}(z, h))^{-1} \frac{d}{dz}{\mathcal K}(z, h) dz.$$
Set $g_h(z): = \det (I + {\mathcal K}(z, h))$ and denote by $M_{\delta}(h)$ the number of the poles $\{\lambda_k\}$ of $R(\lambda)$ such that $h^2\lambda_k$ are in $\omega_{\delta}$. Similarly, we denote by $M_{\delta}^{(j)}(h)$ the number of the eigenvalues $\nu_k$ of $G_D^{(j)}$ such that  $h^2 \nu_k \in \omega_{\delta}.$
Then using the well-known formula
$${\rm tr}\:\Bigl[ (I + {\mathcal K}(x, h))^{-1} \frac{\partial {\mathcal K}(z, h)}{\partial z}\Bigr] = \frac{\partial}{\partial z} \log \det (I + {\mathcal K}(z, h)),$$
we get from (\ref{eq:5.5}) the following
\begin{lemma} Let $\delta \subset Z$ be closed positively oriented curve which avoid the eigenvalues of $h^2 G_D^{(j)},j =1,2$ as well as the poles of $T(z/h^2)^{-1}.$ Then we have
\begin{equation} \label{eq:5.6}
M_{\delta}(h) =   M_{\delta}^{(1)} (h) + M_{\delta}^{(2)}(h) + \frac{1}{2 \pi i} \int_{\delta} \frac{d}{dz} \log g_h(z) dz.
\end{equation}
\end{lemma}

The leading term in (\ref{eq:5.4}) is obtained from the $ M_{\delta}^{(1)} (h) + M_{\delta}^{(2)}(h)$ after a scaling. The crucial point is to examine the asymptotic of the integral involving $\log g_h(z).$ The details of this analysis 
are given in \cite{PV}.


\begin{thebibliography}{99}

\bibitem{AV} M. S. Agranovich and M. I. Vishik, {\em Elliptic problems with
a parameter and parabolic problems of general type}, (Russian)  Uspehi
Mat. Nauk, {\bf  19}, (3) (1964), 53-161. 

\bibitem{BK}  K. Kh. Boimatov and A. G. Kostyuchenko,
{\em Spectral asymptotics of nonselfadjoint elliptic systems of differential
operators in bounded domains}, Matem. Sbornik, {\bf 181}
(1990) (in Russian), English translation: Math. USSR Sbornik {\bf 71} (1992), 517-553.

\bibitem{CH} F. Cakoni and H. Haddar, {\em Transmission eigenvalues in Inverse Scattering Theory}. 
Inside Out II, MSRI Publications {\bf 60}, 2012.

\bibitem{CPV} F. Cardoso, G. Popov and G. Vodev, {\em Asymptotics of the number of resonances in the transmission problem}, 
Commun. Partial Diff. Equations {\bf 26} (2001), 1811-1859.

\bibitem{CPR1} F. Colombini, V. Petkov and J. Rauch, {\em Spectral problems for non elliptic symmetric systems with dissipative boundary conditions}, J. Funct. Anal. {\bf 267} (2014), 1637-1661.  

\bibitem{CPR2} F. Colombini, V. Petkov and J. Rauch, {\em Eigenvalues for  Maxwell's equations with dissipative boundary conditions}, Asymptotic Analysis, {\bf 90} (1-2) (2016), 105-124.

\bibitem{DS} M. Dimassi and J. Sj\"ostrand,  {\em Spectral asymptotics in semi-classical limit}, London Mathematical Society, Lecture Notes Series, {\bf 268}, Cambridge University Press, 1999.

\bibitem{HK} M. Hitrik, K. Krupchyk, P.  Ola, L.  P\"aiv\"arinta, {\em The interior transmission problem and bounds of transmission
eigenvalues}, Math. Res. Lett. {\bf 18}, (2011), 279-293.

\bibitem{LV} E. Lakshtanov and B. Vainberg, {\em Remarks on interior transmission 
eigenvalues, Weyl formula and branching billiards}, J. Phys. A: Math. Theor. {\bf 45} (2012), 125202.

\bibitem{LP1}  P. Lax and R. Phillips, {\em Scattering Theory}, 2nd Edition,  Academic Press, New York, 1989.

\bibitem{LP2} P. Lax and R. Phillips, {\em Scattering theory for dissipative systems}, J. Funct. Anal. {\bf 14} (1973), 172-235.

\bibitem{Ma} A. Majda, {\em The location of the spectrum for the dissipative acoustic operator}, Indiana Univ. Math. J. {\bf 25} (1976), 973-987.

\bibitem{P1} V. Petkov, {\em Location of the eigenvalues of the wave equation with dissipative boundary conditions}, Inverse Problems and Imaging, to appear, (arXiv: math.AP. 1504.06408v4).

\bibitem{PV} V. Petkov and G. Vodev, {\em Asymptotics of the number of the interior transmission eigenvalues}, Journal of Spectral Theory, to appear.

\bibitem{PV1} V. Petkov and G. Vodev, {\em Location of the interior transmission eiegnvalues for a ball}, (arXiv: math.AP. 1603.04604).

\bibitem{PS} H. Pham and P. Stefanov, {\em Weyl asymptotics of the transmission eigenvalues for a constant index of
refraction}, Inverse Problems and Imaging {\bf 8} (3) (2014), 795-810.

\bibitem{R} L. Robbiano, {\em Counting function for interior transmission eigenvalues},  Mathematical Control and Related Fields, 6 (1) (2016), 167-183.

\bibitem{S} J. Sylvester, {\em Transmission eigenvalues in one dimension}, Inverse Problems {\bf 29} (2013), 104009.

\bibitem{SjV} J. Sj\"ostrand and G. Vodev, {\em Asymptotics of the number of Rayleigh resonances}, Math. Ann. {\bf 309} (1997), 287-306.

\bibitem{Sj} J. Sj\"ostrand, {\em Weyl law for semi-classical resonances with randomly perturbed potentials}, M\'emoire de SMF, {\bf 136} (2014).

\bibitem{V1} G. Vodev, {\em Transmission eigenvalue-free regions}, Commun. Math. Phys. {\bf 336} (2015), 1141-1166.

\bibitem{V2} G. Vodev, {\em Transmission eigenvalues for strictly concave domains}, Math. Ann. DOI 10.1007/s00208-015-1329-2, to appear.

\end{thebibliography}
\end{document}